\documentclass{article}
\usepackage{amsmath,amssymb,bbm,theorem,stmaryrd, multicol,verbatim}
\usepackage[french]{babel}
\usepackage[utf8]{inputenc}
\usepackage{hyperref}

\input xy
\xyoption{all}

\def\ZZ {\mathbbm Z}

\newtheorem{thm}{Théorème}

\begin{document}
\title{Calculs effectifs de congruences entre caractères de Dirichlet}
\author{J.~Puydt\\
Institut Joseph Fourier UMR5582\\
Grenoble, France\\
E-mail:\texttt{julien.puydt@ujf-grenoble.fr}
}

\maketitle

\begin{abstract}
  This article aims to find explicit congruences between Dirichlet
  characters and gives various results on how to find some effectively
  on a computer. It ends with concrete examples putting those ideas in
  application.
\end{abstract}

\section{Introduction}
\label{sec:intro}

Le but de ce texte est, à partir du choix de deux entiers
$N\geqslant2$ et $M\geqslant2$, d'obtenir une congruence modulo $m$
pour les caractères de Dirichlet modulo $N$.

Plus précisément, si $G$ est le groupe des caractères de Dirichlet
modulo $N$, alors on cherche des entiers algébriques explicites
$(\alpha_\chi)_{\chi\in G}$ tels que:
\[
\forall x\in\ZZ,
\sum_{\chi\in G}\alpha_\chi\chi(x)
\equiv
0
\bmod M
\]

La motivation derrière ce travail est le théorème suivant, que l'on
trouve dans la thèse~\cite{Puydt2003a} de l'auteur:
\begin{thm}
  On suppose que l'on dispose d'une combinaison linéaire finie de
  caractères de conducteurs $M_S^{\underline\nu}$, qui vérifie une
  congruence:
  \[
  \sum_\chi\gamma_\chi\chi
  \equiv
  0
  \bmod
  M_S^{\underline\tau}\mathcal O
  \]
  alors:
  \[
  C^\alpha_S
  \sum_\chi\gamma_\chi
  H(f,\alpha,S,j,\chi)
  L_{f_{\alpha,S,0}}(j+1,\overline{\chi})
  \equiv
  0
  \bmod
  M_S^{(j-h)\underline\nu+\underline\tau}
  \mathcal O
  \]
où:
\begin{eqnarray*}
  H(f,\alpha,S,j,\chi)
  &=&
  \frac
  {a_{1,0}G_\chi
    \Gamma(k-1)(2i\pi)^{j+1}(-1)^k
    c_\chi^j\alpha_S^{-\underline{\nu}}}
  {\Gamma(j+1)\langle f^0_{\alpha, S}, f_{\alpha, S, 0}\rangle_{N_S}}
  \\
  & &
  \times
  \zeta_{j+1,0}(1_{\hat{\mathbb Z}}\times\mathcal F^{-1}\overline{\chi})(1_f)
  L_{f_{\alpha,S,0}}(k-1,\mathcal F1_{Y_S})
\end{eqnarray*}
\end{thm}

Ce théorème abstrait sur lequel on ne souhaite pas s'étendre (il y a
beaucoup de notations et de concepts qui sont loin du sujet du présent
article), affirme l'existence de certaines congruences pour des
valeurs spéciales de fonctions arithmétiques (des fonctions $L$ de
formes modulaires, tordues par des caractères), qui sont vérifiées si
l'on dispose de congruences pour les caractères. On souhaite donc
construire des exemples concrets de ces dernières, pour un jour
vérifier expérimentalement dans quelle mesure les congruences entre
les valeurs spéciales sont optimales.

Le plan de ce travail est le suivant:
\begin{itemize}
\item on commence par expliquer comment ramener le problème initial à
  un calcul de noyau matriciel ;
\item on explique ensuite comment déterminer le noyau, avec le
  logiciel SAGE\footnote{version 5.4, disponible sur
    \texttt{http//www.sagemath.org}}), d'une façon complètement
  automatique d'une part (mais très inefficace) et de façon plus
  guidée d'autre part (mais très manuelle) ;
\item on détermine des congruences dans des cas particuliers
  explicites ;
\item finalement, on donne le code source des fonctions utilisées.
\end{itemize}

Je souhaite remercier Alexei Pantchichkine pour son soutien
indéfectible, et William Stein pour m'avoir donné accès aux machines
du réseau math.washington.edu\footnote{Les calculs présentés dans cet
  article sont tous aisés et rapides sur une machine simple, mais les
  tâtonnements et expériences nécessaires à leur recherche ont parfois
  nécessité l'accès à une puissance de calcul plus conséquente ; les
  machines en question ont été financées par "National Science
  Foundation Grant No. DMS-0821725".}.

\section{Réduction à un calcul de noyau matriciel}
\label{sec:kernel_matrix}

On se donne comme dans l'introduction un premier entier $N\geqslant2$,
et on considère le groupe $G$ des caractères de Dirichlet modulo $N$,
que l'on liste sous la forme $\chi_1,\dots,\chi_n$ où $n=|G|$. Si on
définit $m=N-1$, on peut alors représenter les valeurs prises par tous
les caractères simultanément sous la forme d'une matrice de taille
$(m,n)$, avec:
\[
\forall(i,j),
A_{i,j}
=
\chi_j(i-1)
\]

De cette façon, si $(\alpha_1,\dots,\alpha_n)$ est une famille de
scalaires et $x\in\llbracket0,N-1\rrbracket$, la quantité
$\sum_{i=1}^n\alpha_i\chi_i(x)$ est le coefficient en ligne $x+1$ dans
le vecteur colonne $A(\alpha_1 \dots \alpha_n)^T$.

Notons que cette matrice initiale $A$ a un noyau trivial: c'est le
résultat bien connu d'indépendance des caractères.

Cependant, comme les caractères considérés sont à valeurs dans les
racines de l'unité, si on fixe $F$ un corps de nombre qui les
contient, la matrice $A$ est à coefficients dans l'anneau des entiers
$\mathcal O$ de ce corps. Pour des raisons pratiques, on choisira bien
sur le corps le plus petit possible, donc un corps cyclotomique de
plus petit degré possible (que l'on notera $d$).

Si on se donne maintenant un entier $M\geqslant2$, la recherche de
congruences modulo cet entier à coefficients dans $\mathcal O$ revient
à chercher un vecteur-colonne $V^T=(\alpha_1 \dots \alpha_n)$ tel que
$AV\in M\mathcal O^n$. Si on note $B$ l'image de $A$ dans l'anneau des
matrices à coefficients dans $\mathcal O/(M)$, on est ramené à
chercher les éléments du noyau de $B$, puis à les relever.

Le problème des congruences pour les caractères de Dirichlet est ainsi
ramené à un problème linéaire de calcul de noyau de matrice à
coefficients dans un quotient d'anneau d'entiers.

\section{Calcul automatique du noyau}
\label{sec:calcul_ker_auto}

D'après un résultat classique (voir par exemple le lemme 6 de
l'article de Birch~\cite{Birch1965a}), l'anneau $\mathcal O$ est le
$\ZZ$-module engendré par les puissances d'une racine primitive de
l'unité. L'anneau $\mathcal O/(M)$ est donc un module de type fini sur
$\ZZ/M\ZZ$, donc un ensemble fini: il suffit de considérer tous les
vecteurs possibles pour obtenir le noyau.

C'est théoriquement très satisfaisant, mais en pratique: $\ZZ/M\ZZ$
est de cardinal $M$, donc $\mathcal O/(M)$ est de cardinal $M^d$, et
les vecteurs à considérer sont donc $M^{dn}$.

Pour le cas le plus simple que l'on présentera de façon détaillée
en~\ref{sec:N5M16}, $N=5$ et $M=16$, cela représente $4294967296$
vecteurs à considérer. Comme de plus les calculs vectoriels dans
l'anneau en question sont en plus relativement lents, il est clair que
cette approche est d'une utilité assez limitée.

Donnons malgré tout le code source effectuant le calcul ; il utilise
une fonction {\texttt matrix\_of\_Dirichlet\_group} discutée plus loin:
\begin{verbatim}
def integer_mod_iter(d,M,zeta):
    for coeffs in product(*tee(range(M),d)):
        result = 0
        for ii in range(d):
            result=result+coeffs[ii]*zeta**ii
        yield result

N = 5
M = 16
A = matrix_of_Dirichlet_group(N)
O = A.base_ring().ring_of_integers()
d = A.base_ring().degree()
Q = O.quotient_ring(O.ideal(M), 'a')
zeta = Q.gens()[1]
B = matrix(Q,A)
m,n = B.dimensions()

print(M**(d*n))
for X in product(*tee(integer_mod_iter(d,M,zeta), n)):
    img = B*vector(Q,X)
    if img.is_zero():
        print(X)
\end{verbatim}

\section{Calcul guidé du noyau}
\label{sec:calcul_ker_guide}

\subsection{Principe général}
\label{sec:principe}

Le problème du calcul du noyau d'une matrice à coefficients dans un
corps est simple ; on peut par exemple penser à un passage sous forme
de matrice échelonnée réduite (c'est, par exemple car on le trouve
partout, l'algorithme~7.3 dans le livre de
Stein~\cite{Stein2007a}). On peut aussi songer à la notion de diviseur
élémentaire, mais elle nécessite un anneau principal et le calcul
effectif un anneau euclidien.

La matrice qui nous intéresse est à coefficients dans un anneau
quotient, qui n'a pas de très bonnes propriétés: les fonctions
habituelles et les algorithmes implémentés dans SAGE ne fonctionnent
donc pas.

L'idée va être d'obtenir une relation de la forme $B*R=L*E$, où $L$ et
$R$ sont des matrices inversibles, et $E$ est aussi proche que
possible d'une matrice échelonnée réduite. Le calcul de vecteurs dans
le noyau de $E$ fournit alors des vecteurs dans le noyau de $B$, tout
simplement en considérant leur image par $R$.

Plus précisément, on va tenter d'obtenir $E$ sous forme de blocs:
\[
\left(\begin{array}{r|r|r}
I_r &
0 &
0 \\
\hline
0 &
Q &
0 \\
\hline
0 &
0 &
0 
\end{array}\right)
\]

où le bloc $I_r$ est un bloc identité (donc carré), et $Q$ est un bloc
(a priori rectangulaire) dans lequel on garantit:
\begin{itemize}
\item aucun coefficient n'est inversible ;
\item aucune ligne n'est nulle ;
\item aucune colonne n'est nulle.
\end{itemize}

C'est le mieux que l'on puisse demander à SAGE de faire
automatiquement ; c'est la raison pour laquelle la suite de la
recherche nécessite une intervention manuelle (discutée plus loin).

Le calcul du noyau est alors ramené au calcul du noyau de cette
matrice réduite $E$ ; il fait intervenir:
\begin{itemize}
\item des éléments triviaux, qui correspondent aux colonnes nulles à
  droite (s'il y en a) ;
\item des éléments non triviaux, qui correspondent au noyau du bloc
  $Q$ (que l'on voudra donc réduire à la main après avoir obtenu une
  première version).
\end{itemize}

On appelera l'entier $r$ (taille de la matrice identité dans la
décomposition) le \emph{pseudo-rang}: il correspond à un nombre de
colonnes que l'on sait être indépendantes sur l'anneau considéré. Il
dépend évidemment de la décomposition obtenue, mais s'il est égal à
$n$, on est assuré qu'il n'existe pas de congruences modulo $M$. Le
nombre de colonnes nulles sera appelé le \emph{noyau garanti}.

\subsection{Algorithme de réduction automatique}
\label{sec:algo_reduction}

L'idée est d'aménager l'algorithme du pivot de Gauss, en l'analysant
en termes d'opérations élémentaires sur les lignes et les colonnes.

Si $(i,j)$ est un couple d'indices distincts, on peut définir
$P_{i,j}$ la matrice de permutation associé à la transposition de $i$
et $j$, et si $a$ est un scalaire, on peut considérer $T_{i,j}(a)$ la
matrice de transvection qui est une matrice carrée identité modifiée
avec le coefficient $a$ en ligne $i$ et colonne $j$, et si $a\neq0$,
la matrice de dilatation $D_i(a)$, matrice identité modifiée avec le
coefficient $a$ en ligne $i$ et colonne $i$. Les opérations
élémentaires sur les lignes se décrivent par multiplication à gauche
par ces matrices, et les opérations élémentaires sur les colonnes par
multiplication à droite. Attention, ces matrices sont carrées, mais on
ne fait pas apparaître leur taille dans la notation.

Le principe va être le suivant: si on a réussi à écrire $B*R=L*E$ avec
$E$ quelconque, alors on a aussi les égalités suivantes, pour les
opérations sur les lignes:
\begin{eqnarray*}
  B*R
  & = &
  LP_{i,j}*P_{i,j}E
  \\
  & = &
  LT_{i,j}(-a)*T_{i,j}(a)E
  \\
  & = &
  LD_i(1/a)*D_i(a)E
\end{eqnarray*}
et sur les colonnes:
\begin{eqnarray*}
  B*RP_{i,j}
  & = &
  L*EP_{i,j}
  \\
  B*RT_{i,j}(a)
  & = &
  L*ET_{i,j}(a)
  \\
  B*RD_i(a)
  & = &
  L*ED_i(a)
\end{eqnarray*}
qui montrent que si une opération élémentaire permet de \og simplifier
$E$\fg{}, on peut le faire quitte à modifier $L$ (pour les opérations
sur les lignes, et par une opération sur les colonnes) ou $R$ (pour
les opérations sur les colonnes, et par la même opération sur les
colonnes).

Expliquons comment on travaille:
\begin{enumerate}
\item On part tout d'abord avec le triplet $(L,E,R)=(I_m,B,I_n)$, qui
  vérifie par construction $B*I_n=I_m*E$ ; c'est l'invariant $B*R=L*E$
  de l'algorithme. On commence avec $k=1$ ;
\item Si on trouve un coefficient inversible dans $E$ en position
  $(l,c)$, alors:
  \begin{enumerate}
  \item on permute les lignes $l$ et $k$ de $E$, et les colonnes $l$
    et $k$ de $L$ (ce qui conserve l'invariant) ;
  \item on permute les colonnes $c$ et $k$ de $E$ de $R$ (pour
    l'invariance) ;
  \item le coefficient unité se retrouve en $(k,k)$ ; par dilatation,
    on divise la ligne $k$ de $E$ et on multiplie la colonne $k$ de
    $L$ par ce coefficient (même remarque) ;
  \item le coefficient est maintenant égal à $1$ ; on ajoute dans $E$,
    avec multiplication par un coefficient la $k$-ème ligne puis la
    $k$-ème colonne aux autres lignes et colonnes pour que la ligne
    $k$ et la colonne $k$ ne contiennent que ce $1$ (c'est un pivot) ;
    opérations sur $E$ que l'on compense sur $L$ et $R$ comme discuté
    ci-dessus, toujours pour conserver l'invariant ;
  \item on incrémente $k$ et on recommence (dans la mesure où
    $k\leqslant\min\{m,n\}$).
  \end{enumerate}
\item Si on ne trouve plus de coefficient inversible, on se contente
  de faire migrer par des permutations les lignes et les colonnes
  nulles de $B$ vers le bas et la droite (toujours en compensant sur
  $L$ et $R$).
\end{enumerate}

La première partie de pivot de l'algorithme garantit que l'on a une
allure:
\[
\left(\begin{array}{r|r}
I_r &
0 \\
\hline
0 &
* 
\end{array}\right)
\]
sans élément inversible hors du premier bloc, et la migration finale
des lignes et colonnes nulles garantit que le bloc central n'a plus ni
ligne ni colonne nulle et fait apparaître les parties droites et
basses nulles.

\subsection{Discussion de l'implémentation}
\label{sec:implementation}

\subsubsection{\textrm matrix\_of\_Dirichlet\_group}

Cette fonction reçoit l'entier $N\geqslant2$ choisi, puis calcule et
renvoie la matrice associée au groupe des caractères de Dirichlet,
qu'il faudra ensuite réduire modulo un entier $M\geqslant2$. C'est la
seule qui ne soit pas contenue dans l'implémentation de la classe.

\subsubsection{\textrm dirty\_cached\_is\_unit}

La version du logiciel SAGE utilisée ne dispose pas d'un moyen de test
de l'inversibilité dans un anneau quotient de l'anneau des entiers
d'un corps de nombres ; mais propose une méthode d'inversion sous
réserve d'existence (avec levée d'exception si ce n'est pas le
cas). On a donc implémenté une fonction de test à l'aide de
l'inversion.

C'est très coûteux, mais aura le mérite d'être aisé à remplacer lors
d'évolutions futures du logiciel. Pour gagner un peu de temps, on
mémorise les résultats. Expéri\-mentalement, on gagne relativement peu
avec cette mise en cache (quelques pourcents).

Elle retourne un booléen.

\subsubsection{Classe {\textrm invariant\_triplet}}

On travaille comme on l'a vu avec des triplets de matrices
$(L,E,R)$. La présence de $L$ n'a aucun intérêt a priori vu l'objectif
de calcul d'un noyau, mais on la conserve tout de même, pour deux
raisons:
\begin{itemize}
\item la première est qu'elle coûte peu à maintenir donc ne constitue
  pas une réelle gêne,
\item et la seconde parce que le maintien de l'invariant donne un
  moyen de vérification intéressant. En particulier, durant le
  développement, de nombreuses coquilles ont été repérées par des
  tests de l'invariant sur des matrices tirées au hasard.
\end{itemize}

Ensuite, comme on le verra plus loin, l'algorithme automatique ne
donne pas de très bons résultats: une fois choisi le couple $(M,N)$,
il faut guider la machine. C'est la raison pour laquelle on a modélisé
la situation à l'aide d'une classe {\textrm invariant\_triplet}, dont
on initialise les instances avec la matrice $B$, et qui va se charger
de gérer le triplet.

Détaillons l'interface programmatique de cette classe:
\begin{description}

\item[Initialisation] La fonction d'initialisation reçoit une matrice
  de taille quelconque, et initialise l'objet avec le triplet
  $(L,E,R)=(I_m,B,I_n)$.

\item[Variables] On peut inspecter les variables $L$, $R$ et $E$. En
  particulier, l'inspection de $E$ permet de décider quelles
  opérations de réduction on va demander.

\item[Méthode {\textrm assert\_invariant}] Vérifie que l'objet
  respecte toujours l'invariant -- extrèmement utile durant le
  développement.

\item[Méthode {\textrm check\_kernel\_vector}] Reçoit une liste de
  coefficients, qu'elle utilise pour créer un vecteur-colonne. Si ce
  vecteur colonne est dans le noyau de $E$, alors la méthode retourne
  une paire constituée du booléen {\textrm True} et de l'image du
  vecteur-colonne par la matrice $R$ (donc le vecteur {\bf explicite}
  donnant la congruence cherchée!) ; sinon elle retourne une paire
  constituée du booléen {\textrm False} et de l'image du
  vecteur-colonne par la matrice $E$ (ce qui permet de voir à quel
  point on s'est trompé pour corriger).

\item[Méthodes {\textrm do\_row\_addition} et {\textrm
    do\_column\_addition}] reçoivent deux indices $i$ et $j$ et un
  coefficient $a$. Elles réalisent l'action de $T_{i,j}(a)$ sur le
  triplet.

\item[Méthodes {\textrm do\_swap\_rows} et {\textrm
    do\_swap\_columns}] recoivent deux indices $i$ et $j$, et
  réalisent l'action de $P_{i,j}$ sur le triplet.

\item[Méthode {\textrm do\_pivot}] réalise la partie pivot de l'algorithme de
  réduction présenté en~\ref{sec:algo_reduction},
  page~\pageref{sec:algo_reduction}.

\item[Méthode {\textrm do\_migrate\_zeros}] réalise la partie de déplacement des
  lignes et des colonnes nulles vers le bas et la droite dans
  l'algorithme présenté en~\ref{sec:algo_reduction}.

\item[Méthode {\textrm do\_normalize}] appelle successivement les fonctions
  {\textrm do\_pivot} et {\textrm do\_migrate\_zeros}.

\end{description}

\section{Résultats obtenus}
\label{sec:resultats}

\subsection{Statistiques}
\label{sec:stats}

La réduction automatique, passant en revue les cas où
$M\in\llbracket2,20\rrbracket$ et $N\in\llbracket2,20\rrbracket$ (soit
$361$ paires), a trouvé $64$ cas où le noyau est garanti comme réduit
à zéro (matrice $I_r$ occupant toute la largeur de la matrice), et où
donc aucune congruence n'est possible.

Les fréquences des rangs garantis parmi les $297$ cas restants est
listée dans la table~\ref{freq:rangs} et les fréquences des noyaux
garantis dans la table~\ref{freq:noyaux}. On constate donc que les
matrices rencontrées même avec de petits paramètres ne se réduisent
pas bien automatiquement: la recherche de congruences nécessite
beaucoup d'affinage manuel.

\begin{table}[h]
\caption{\label{freq:rangs} Fréquences des pseudo-rangs}
\centering
\begin{tabular}{|c|c|}
  \hline
  pseudo-rang & fréquence \\
  \hline
  1 & 180 \\
  2 & 54 \\
  3 & 27 \\
  6 & 27 \\
  7 & 9 \\
  \hline
\end{tabular}
\end{table}

\begin{table}[h]
\caption{\label{freq:noyaux} Fréquences de noyaux}
\centering
\begin{tabular}{|c|c|}
  \hline
  noyau garanti & fréquence \\
  \hline
  0 & 279 \\
  1 & 13 \\
  3 & 5 \\
  \hline
\end{tabular}
\end{table}

\subsection{Congruences pour $N=5$ et $M=16$}
\label{sec:N5M16}

On choisit de se concentrer sur le cas $N=5$ et $M=16$ car c'est un
cas de congruence de degré $4$ le long de $2$, et la dimension est
petite: cela semble donc un bon choix de premier exemple. Remarquons
que cet exemple minimaliste pour la méthode guidée est hors de portée
de la recherche exhaustive, comme on l'a vu
en~\ref{sec:calcul_ker_auto}.

Voici le début de la session de calcul:
\begin{verbatim}
N = 5
M = 16
A = matrix_of_Dirichlet_group(N)
O = A.base_ring().ring_of_integers()
m,n = A.dimensions()
Q = O.quotient_ring(O.ideal(M), 'a')
B = invariant_triplet(matrix(Q, A))
B.do_normalize()
\end{verbatim}

Que l'on poursuit ensuite à la main par:
\begin{verbatim}
zeta4=Q.gens()[1]
B.do_row_addition(3,2,1)
B.do_row_addition(3,1,-1)
B.do_column_addition(3,1,1)
B.do_row_addition(1,2,-2)
B.do_row_addition(1,3,1)
\end{verbatim}

À partir de là, la matrice $E$ est assez simplifiée pour que l'on
devine des vecteurs du noyau ; après avoir éliminé les doublons, il
reste les vecteurs:
\begin{verbatim}
B.check_kernel_vector([0,0,0,8])
B.check_kernel_vector([0,0,4,0])
B.check_kernel_vector([0,8,0,4*zeta4-4])
\end{verbatim}

On trouve donc des vecteurs de congruences explicites:
\[
\left(\begin{array}{c}
0 \\ 8 \\ 0 \\ 8
\end{array}\right)
,
\left(\begin{array}{c}
4 \\ -4 \\ 4 \\ -4
\end{array}\right)
,
\left(\begin{array}{c}
8\zeta_4 \\ 4\zeta_4+4 \\ 0 \\ 4\zeta_4-4
\end{array}\right)
\]

\subsection{Congruences pour $N=7$ et $M=15$}
\label{sec:N7M15}

Voici le début de la session de calcul:
\begin{verbatim}
N = 7
M = 15
A = matrix_of_Dirichlet_group(N)
O = A.base_ring().ring_of_integers()
m,n = A.dimensions()
Q = O.quotient_ring(O.ideal(M), 'a')
B = invariant_triplet(matrix(Q, A))
B.do_normalize()
\end{verbatim}

(a parte: $M^{dn}$ vaut ici $129746337890625$)


Que l'on poursuit ensuite à la main ; au bout de quelques étapes, on
est ramené à:
\[
E=
\begin{pmatrix}
  1 & 0 & 0 & 0 & 0 & 0 \\
  0 & 1 & 0 & 0 & 0 & 0 \\
  0 & 0 & 1-2\zeta6 & 0 & 2-\zeta_6 & 0 \\
  0 & 0 & 0 & 2(2\zeta_6-1) & 0 & -2(\zeta_6+1) \\
  0 & 0 & 0 & 0 & 3 & 0 \\
  0 & 0 & 0 & 0 & 0 & -6 \\
  0 & 0 & 0 & 0 & 0 & 0 \\
\end{pmatrix}
\]
et:
\[
R=
\begin{pmatrix}
  1 & -1 & 0 & 0 & 0 & -1 \\
  0 & 1 & -1 & 1 & -1 & -1 \\
  0 & 0 & 0 & 1 & 0 & 0 \\
  0 & 0 & 0 & 0 & 1 & 1 \\
  0 & 0 & 0 & -1 & 0 & 1 \\
  0 & 0 & 1 & -1 & 0 & 0
\end{pmatrix}
\]

On trouve donc des vecteurs de congruences explicites:
\[
\left(\begin{array}{c}
0 \\ 5(\zeta_6+1) \\ 5(\zeta_6+1) \\ 0 \\ -5(\zeta_6+1) \\ -5(\zeta_6+1)
\end{array}\right)
\left(\begin{array}{c}
0 \\ -5(\zeta_6+1) \\ 0 \\ 0 \\ 0 \\ 5(\zeta_6+1)
\end{array}\right)
,
\left(\begin{array}{c}
0 \\ 5 \\ 0 \\ 5 \\ 0 \\ 5
\end{array}\right)
,
\left(\begin{array}{c}
-5 \\ 5 \\ -5 \\ 5 \\ -5 \\ 5
\end{array}\right)
\]

\section{Code}
\label{sec:code}

\begin{verbatim}
def matrix_of_Dirichlet_group (M):
    G=DirichletGroup(M)
    characters = G.list()
    field = G.base_ring()
    return matrix(field, M, len(characters), lambda lig,col: characters[col](lig))

cached_units = {}
def dirty_cached_is_unit(a):
    if cached_units.has_key(a):
        return cached_units[a]
    else:
        result = False
        O = a.base_ring()
        try:
            b = O(1/a)
            result = True
        except:
            pass
        cached_units[a]=result
        return result

class invariant_triplet:

    def __init__(self, B):
        self.B = B
        self.E = copy(B)
        self.m, self.n = B.dimensions()
        self.L = identity_matrix(B.base_ring(), self.m)
        self.R = identity_matrix(B.base_ring(), self.n)
        self.assert_invariant()

    def assert_invariant(self):
        if self.B*self.R != self.L*self.E:
            raise ArithmeticError

    def check_kernel_vector(self, coeffs):
        v = matrix(self.E.base_ring(), self.n, 1, coeffs)
        Ev = self.E*v
        if Ev.is_zero():
            return (True, self.R*v)
        else:
            return (False, Ev)

    def do_row_addition(self,i,j,a):
        self.E.add_multiple_of_row(i,j,a)
        self.L.add_multiple_of_column(j,i,-a)

    def do_column_addition(self,i,j,a):
        self.E.add_multiple_of_column(i,j,a)
        self.R.add_multiple_of_column(i,j,a)

    def do_swap_rows(self,i,j):
        self.E.swap_rows(i,j)
        self.L.swap_columns(i,j)

    def do_swap_columns(self,i,j):
        self.E.swap_columns(i,j)
        self.R.swap_columns(i,j)

    def do_pivot(self):
        k = 0
        while k < min(self.m,self.n):
            found_pivot = False
            l = k
            c = k
            while l < self.m and c < self.n and not found_pivot:
                if dirty_cached_is_unit(self.E[l,c]):
                    found_pivot = True
                if not found_pivot:
                    l = l+1
                    if l == self.m:
                        l = k
                        c = c+1
            if found_pivot:
                self.do_swap_rows(l,k)
                self.do_swap_columns(c,k)
                coeff=self.E[k,k]
                self.E.rescale_row(k,1/coeff)
                self.L.rescale_col(k,coeff)
                for l in range(self.m):
                    if l != k and self.E[l,k] != 0:
                        coeff=self.E[l,k]
                        self.do_row_addition(l,k,-coeff)
                for c in range(self.n):
                    if c != k and self.E[k,c] != 0:
                        coeff = self.E[k,c]
                        self.do_column_addition(c,k,-coeff)
                k = k+1
            else:
                break

    def do_migrate_zeros(self):
        target = self.m-1
        considered = target
        while considered >= 0:
            if self.E[considered,:].is_zero():
                if considered < target:
                    self.E.swap_rows(considered, target)
                    self.L.swap_columns(considered, target)
                target = target-1
            considered = considered-1

        target = self.n-1
        considered = target
        while considered >= 0:
            if self.E[:,considered].is_zero():
                if considered < target:
                    self.E.swap_columns(considered, target)
                    self.R.swap_columns(considered, target)
                target = target-1
            considered = considered-1

    def do_normalize(self):
        self.do_pivot()
        self.do_migrate_zeros()
\end{verbatim}

\bibliography{../Abords/bibliography}

\begin{thebibliography}{1}
\expandafter\ifx\csname fonteauteurs\endcsname\relax
\def\fonteauteurs{\scshape}\fi

\bibitem{Birch1965a}
B.~J. \bgroup\fonteauteurs\bgroup Birch\egroup\egroup{} :
\newblock {C}yclotomic fields and {K}ummer extensions.
\newblock \emph{In} J.~W.~S. \bgroup\fonteauteurs\bgroup
  Cassels\egroup\egroup{} et A.~\bgroup\fonteauteurs\bgroup
  Fröhlich\egroup\egroup{}, \'editeurs :  {\em {A}lgebraic number theory},
  pages 85--93. {A}cademic press, 1965.

\bibitem{Puydt2003a}
J.~\bgroup\fonteauteurs\bgroup Puydt\egroup\egroup{} :
\newblock {\em {V}aleurs spéciales de fonctions {$L$} de formes modulaires
  adéliques}.
\newblock Th\`ese de doctorat, {I}nstitut {F}ourier, 2003.

\bibitem{Stein2007a}
W.~A. \bgroup\fonteauteurs\bgroup Stein\egroup\egroup{} :
\newblock {\em {M}odular forms, a computational approach}, volume~79 de {\em
  {G}raduate studies in mathematics}.
\newblock {A}merican mathematical society, 2007.

\end{thebibliography}
\bibliographystyle{plain-fr}

\end{document}